\documentclass [12pt]{article}      
\usepackage{amsmath} 
\usepackage{amsthm}
\usepackage{latexsym}
\usepackage{pb-diagram}
\usepackage{amssymb}
\usepackage{bbm}
\usepackage{color}
\usepackage{graphicx}
\usepackage{hyperref}
\hypersetup{
    colorlinks,
    citecolor=blue,
    filecolor=blue,
    linkcolor=blue,
    urlcolor=blue
}
\hypersetup{linktocpage}
\usepackage{enumerate}
\usepackage{todonotes}
\makeatletter
\let\r@eqnnum\@eqnnum
\input{leqno.clo}
\let\l@eqnnum\@eqnnum
\newcommand{\leqnos}{\let\@eqnnum\l@eqnnum}
\newcommand{\reqnos}{\let\@eqnnum\r@eqnnum}
\reqnos
\makeatother

\newtheorem{theorem}{Theorem}
\newtheorem*{theorem*}{Theorem}
\newtheorem{corollary}[theorem]{Corollary}
\newtheorem*{corollary*}{Corollary}

\newtheorem{lemma}[theorem]{Lemma}
\newtheorem{definition}[theorem]{Definition}
\newtheorem{prop}[theorem]{Proposition}
\newtheorem{proposition}[theorem]{Proposition}

\newtheorem{example}[theorem]{Example}

\newtheorem{remark}[theorem]{Remark}

\newcommand{\cat}{^\frown}
\newcommand{\rest}{\ensuremath{\upharpoonright}}

\newcommand{\inv}{{^{-1}}}

\newcommand{\la}{\langle}
\newcommand{\ra}{\rangle}

\newcommand{\Chi}{\ensuremath{\mbox{\Large{$\chi$}}}}

\newcommand{\xbmt}{\ensuremath{(X,\mcb,\mu,T)}}

\newcommand{\sz}{\Sigma^\poZ}

 \newcommand{\poN}{\mathbb N}
\newcommand{\poZ}{\mathbb Z}
\newcommand{\nn}{{\mathbb N}}

\newcommand{\poz}{\mathbb Z}
\newcommand{\bk}{{\mathbb K}}
\newcommand{\bl}{\mathbb L}

\newcommand{\mct}{{\ensuremath{\mathcal T}}}

\newcommand{\mcq}{\ensuremath{\mathcal Q}}

\newcommand{\mcp}{\mathcal P}

\newcommand{\mcw}{\mathcal W}

\newcommand{\mcc}{{\mathcal C}}

\newcommand{\mcb}{\mathcal B}

\newcommand{\mco}{\mathcal O}

\newcommand{\mcv}{\mathcal V}

\newcommand{\bfni}[1]{\noindent {{\bf{#1}}}}

\renewcommand{\qed}{{\nopagebreak \hfill $\dashv$ 
 \par\bigskip}}

\newcommand{\pf}{{\par\noindent{$\vdash$\ \ \ }}}

\newcommand{\skn}{\la k_n:n\in\nn\ra}

 \title{Odometer Based Systems}
 \author{Matthew Foreman and Benjamin Weiss}
 
\begin{document}
 \maketitle
 \begin{abstract}
  Construction sequences are a general method of building symbolic shifts that
capture cut-and-stack constructions and are general enough to give symbolic
representations of Anosov-Katok diffeomorphisms. We show here that any
finite entropy system that has an odometer factor can be represented as a special class of construction sequences, the odometer based construction
sequences which correspond to those cut-and-stack constructions that do not
use spacers. We also show that any additional property called the ``small
word condition" can also be satisfied in a uniform way.
 \end{abstract}
 
%

\section{Introduction}
  
Construction sequences are a general method of building symbolic shifts that capture cut-and-stack constructions and are general enough to give symbolic representations of Anosov-Katok diffeomorphisms.  This paper studies a special class of construction sequences, the \emph{odometer based construction sequences} that corresponds to those cut-and-stack constructions that don't use spacers.
 
 In \cite{part2} we show that there is a functorial isomorphism between the symbolic systems that are limits of odometer based construction sequences and symbolic systems that are limits of a class of construction sequences called \emph{circular systems}. The circular systems, in turn, can be realized as diffeomorphisms of the 2-torus. As a corollary the qualitative ergodic theoretic structure of the odometer based systems is reflected in the diffeomorphisms of the 2-torus. For example one deduces that there are measure-distal diffeomorphisms of the torus of all countable ordinal heights \cite{partV} and for all Choquet simplices $\mathcal K$, there is a Lebesgue measure preserving ergodic diffeomorphism of the torus that has $\mathcal K$ as its simplex of invariant measures.

To use the functor defined in \cite{part2} one needs to see that the class of transformations isomorphic to limits of odometer based construction sequences is quite rich and complicated. This is the point of the current paper. 
 
 It is a classical theorem of Krieger (\cite{Krieger})
     that an ergodic system with finite entropy has a finite generating partition. This gives a symbolic representation for any such 
     system  and shows that the theory of finite entropy ergodic measure preserving systems coincides with the theory of 
     finite valued ergodic stationary processes $\{X_n\}$. When studying stationary processes $\{X_n\}$ it is often useful to have a block structure, namely
 a way of dividing the indices into a hierarchy of blocks of lengths $k_1, k_1k_2, k_1k_2k_3, \dots$
     in a unique fashion. If this is possible then the process will have as a factor the odometer transformation corresponding to the sequence $\{k_n\}$.
     Our main theorem is that it is always possible to find such a symbolic representation with a rather simple form whenever this necessary condition is
     satisfied. 
 
\begin{theorem*}(See \ref{first presentation} in Section \ref{OBSS}
  Let $\xbmt$ be a measure preserving system with finite entropy. Then $X$ has an odometer factor if and only if $X$ is isomorphic to an odometer based symbolic system.
  \end{theorem*}

The class of ergodic transformations containing an odometer factor is easily characterized spectrally as those transformations whose associated unitary operator has infinitely many eigenvalues of finite multiplicative order. 

 The periodic factors of an ergodic system are the obstructions to the ergodicity of powers of $T$. If $T$ is totally ergodic, i.e. all powers are ergodic,
      then the product of $T$ with any odometer is ergodic. In general we have the following proposition which illustrates the ubiquity of ergodic transformations with odometer factors:

\begin{prop}\label{accidental products}
Given any  ergodic transformation ${\mathbb X}=\xbmt$ either:
\begin{enumerate}

\item $\mathbb X$ has an odometer factor 
\begin{center}
or
\end{center}

\item there is an odometer 
$\mathfrak O$ such that $\mathbb X\times \mathfrak O$ is ergodic (and  $\mathbb X\times \mathfrak O$  has finite entropy if $\mathbb X$ does). 

\end{enumerate}
\end{prop}
\noindent In particular, every finite entropy transformation is a factor of a finite entropy odometer based symbolic system and the finite entropy transformations that have an odometer factor are closed under finite entropy extensions.

We should point out that special symbolic processes with a block structure, called Toeplitz systems, have been well studied 
       from the point of view of topological dynamics.  Downarowicz and Lacroix (\cite{downer2}, theorem 8) showed that every transformation satisfying the hypothesis of our main theorem can be represented as the orbit closure of a Toeplitz sequence. Proposition \ref{after downer}, presents orbit closures of Toeplitz sequences as limits of odometer based construction sequences, giving an alternate proof of our main theorem. The authors were unaware the results in \cite{downer2} when the we obtained the results in this paper.\footnote{In addition the proof offered in \cite{downer2} makes reference for a key result
to \cite{Weiss_reference} in which only a sketch of a more general theorem is given and the specific result they need is not even mentioned there.} 
              
       We also note work of Williams presenting the odometer itself as a limit of a construction  sequence (see Williams, \cite{williams}) as well as the recent work of Adams, Ferenczi, and Petersen \cite{AFP}, which realizes generalized odometers and indeed all rank one systems as ``constructive symbolic rank one systems",
in the terminology of \cite{Feren}.

\paragraph{The structure of this paper} Section \ref{preliminaries} has the basic definitions used in the paper as well as properties of Odometer systems that that we use in the construction. Section \ref{OBSS} contains the proof of our main theorem, Theorem \ref{first presentation}. It begins by pointing out a known fact that odometers cannot be represented topologically as symbolic shifts, in contrast to Theorem \ref{first presentation}, which is in the measure category. As a precursor it then presents the odometer as an odometer based system, describes the plan of the proof and finally gives the proof in detail. 

In Section \ref{Toeplitz Systems} we discuss the connections with Toeplitz systems, showing how to augment a Toeplitz system to get an odometer based system while preserving the simplex of invariant measures. It then follows from a remarkable theorem of Downarowicz \cite{downer} (generalizing work of Williams \cite{williams}) saying that arbitrary simplices of invariant measures can be realized on Toeplitz sequences to see that arbitrary simplices of invariant measures can be realized on limits of odometer construction sequences. 

The applications of this paper require that the odometer based construction sequences in the domain of the isomorphism functor has the frequencies of words decreasing arbitrarily fast.  We call this the \emph{small word property}. In Section \ref{SWP} we define the small word property and show that we can realized odometer based systems continuously in a sequence of small word requirements.

\section{Preliminaries}\label{preliminaries}

An \emph{alphabet} is a finite collection of symbols. A \emph{word} in $\Sigma$ is a finite sequence of elements of $\Sigma$. If $w\in \Sigma^{<\nn}$ is a word, we denote its length by $|w|$.  By $\Sigma^{\poZ}$ we mean doubly infinite sequences of letters in $\Sigma$. This has a natural product topology induced by the discrete topology on $\Sigma$. This topology is compact if $\Sigma$ is finite.  For this paper a \emph{symbolic system} is a closed, shift-invariant $\bk\subseteq \Sigma^\poZ$. 

A collection of words $\mcw$ is  \emph{uniquely readable} if and only if whenever $u, v, w\in \mcw$ and $uv=pws$ then either $p$ or $s$ is the empty word.

We note that we can view both words and elements of $\Sigma^\poZ$ as functions. If $f:A\to B$ and $A'\subseteq A$, the restriction of $f$ to $A'$ is denoted $f\rest A'$.

\subsection{Partitions and Symbolic Systems}\label{parts and symbols}
Let $(X,\mcb,\mu)$ be a standard measure space. An ordered \emph{partition} of $X$ is a set $\mathcal P=\la A_0, A_1,\dots \ra$ such that each $A_i\in\mcb$, 
$A_i\cap A_j=\emptyset$ if $i\ne j$, and $X=\bigcup_i A_i$. We allow our partitions to be finite or countable and identify two partitions $\mcp=\la A_i\ra$, $\mcq=\la B_j\ra$ if for all $i, \mu(A_i\Delta B_i)=0$.

We will frequently refer to  ordered countable measurable partitions simply as \emph{partitions}. A partition is finite iff for all large enough $n, \mu(P_n)=0$.  If $\mcp$ and $\mcq$ are partitions then $\mcq$ \emph{refines} $\mcp$ iff the atoms of $\mcq$ can be grouped into sets $\la S_n:n\in \nn\ra$ such that 
\[\sum_n\mu(P_n\Delta (\bigcup_{i\in S_n}Q_i))=0.\]
In this case we will write that $\mcq\ll\mcp$. A \emph{a decreasing sequence of partitions} is a sequence $\la \mcp_n:n\in\nn\ra$ such that for all $m<n, \mcp_n\ll \mcp_m$. If $A\in \mcb$ is a measurable set and $\mcp$ is a partition then we let $\mcp\rest A$  be the partition of $A$ defined as $\la P_n\cap A:n\in \nn\ra$.

\begin{definition}\label{definition of generation}
Let $(X,\mcb, \mu)$ be a measure space. We will say that a sequence of partitions $\la \mcp_n:n\in\nn\ra$ \emph{generates} (or generates $\mcb$) iff the smallest $\sigma$-algebra containing $\bigcup_n\mcp_n$ is $\mcb$ (modulo measure zero sets). If $T$ is a measure preserving transformation we will write $T\mcp$ for the partition $\la Ta:a\in \mcp\ra$. In the context of a measure preserving $T:X\to X$ we will say that a partition $\mcp$ is a \emph{generator} for $T$ iff $\la T^i\mcp:i\in \poZ\ra$ generates $\mcb$.

\end{definition}

Given a measure preserving system $\xbmt$ and a partition $\mcp$ of $X$, define a map $\phi:X\to \mcp^\poZ$ by setting (for each $a\in \mcp$):
\[\phi(x)(n)=a \mbox{ if and only if }T^nx\in a.\]
The bi-infinite sequence $\phi(x)$ will be called the $\mcp$-name of $x$. The closure of $\phi(X)\subseteq \mcp^\poZ$ is a symbolic system.

Define a measure on $\mcp^\poZ$ by setting $\phi^*(\mu)(A)=\mu(\phi^{-1}[A])$. This is a Borel measure on the symbolic shift $\mcp^\poZ$ and makes $(\mcp^\poZ,\mcc,\nu,sh)$ into a factor of $\xbmt$ (where $\nu=\phi^*(\mu)$).  This factor map is an isomorphism if and only if $\mcb$ is the smallest shift-invariant  $\sigma$-algebra containing all of the sets in $\mcp$ (up to sets of measure zero); i.e. $\mcp$ is a \emph{generator} for $T$. In general the support of $\nu$ is the closure of $\phi(X)$.

\begin{remark}\label{joining factors} Let $\mcp, \mcq$ be partitions of 
$X$. Then $\mcp$ and $\mcq$ determine  factors $Y_\mcp$ and $Y_\mcq$.  
Define $\phi:X\to Y_\mcp\times Y_\mcq$ by setting $\phi(x)=(s_p,s_q)$ where $s_{p}$ is the $\mcp$-name of $x$ and $s_{q}$ is the $\mcq$-name of $x$. Let $\eta=\phi^{*}(\mu)$. Then $(Y_\mcp\times Y_\mcq,\mcc,\eta,sh)$ is isomorphic to the smallest factor of $X$ containing both $Y_\mcp$ and $Y_\mcq$ as factors.
\end{remark}

\subsection{Basic Facts About Odometers} Let
$\la k_i:i\in\poN\ra$ be an infinite sequence of integers with $k_i\ge 2$. Then the sequence
$k_i$ determines an \emph{odometer} transformation with domain the
compact space\footnote{We write $\poZ/\poZ_k$ as $\poZ_k$.}
\[\mbox{{\boldmath$O$}}=_{def}\prod_i\poz_{k_i}.\] 

The space {\boldmath$O$} is naturally
a monothetic compact abelian group,  with the operation of addition and ``carrying
right". We will denote the group element $(1, 0, 0, 0,\dots)$ by $\bar{1}$, and the result of adding
$\bar{1}$ to itself $j$ times by $\bar{j}$.

 The Haar measure on this group can be defined explicitly. Define
a  measure $\nu_i$ on each $\poz_{k_i}$ that gives each point measure
$1/k_i$. Then Haar measure $\mu$ is  the
product measure of the $\nu_i$. 

 The odometer transformation
$\mathcal O:\mbox{{\boldmath $O$}}\to \mbox{{\boldmath $O$}}$ is defined by taking
an $x\in
\prod_i\poz_{k_i}$ and adding the group element $\bar{1}$,
More explicitly, $\mathcal O(x)(0)=x(0)+1($mod $ k_0)$ and $\mathcal
O(x)(1)=x(1)$ unless $x(0)=k_0-1$, in which case we ``carry one" and set
$\mathcal O(x)(1)=x(1)+1($mod $ k_1)$, etc.

The map $\mathcal O:\mbox{{\boldmath $O$}}\to \mbox{{\boldmath $O$}}$ is a topologically minimal,
uniquely ergodic,
invertible homeomorphism that preserves the measure
$\mu$.  When we are viewing the odometer as a measure preserving system we will denote it by 
$\mathfrak O$.

Define $U_{\mathfrak O}:L^2(\mathfrak O)\to L^2(\mathfrak O)$ by setting
$U_{\mathfrak O}(f)=f\circ \mathcal O$. Then $U_{\mathfrak O}$ is the canonical unitary operator associated with $\mathcal O$.
The characters $\Chi\in\hat{O}$ are eigenfunctions for the $U_{\mathfrak O}$ since
\[\Chi(x+\bar{1})=\Chi(\bar{1})\Chi(x).\]
 Since the characters form a basis for $L^2(\mathfrak O)$, the
odometer map has discrete spectrum. 

Here is an explicit description of the characters.
Fix $n$ and let $K_n=\prod_{i<n}k_i$.  Let $A_0\subset \prod_i\poz_{k_i}$ 
be the collection of points
whose first $n+1$ coordinates are zero, and for $0\le k<K_n$ set
$A_k=\mathcal O^k(A)$. Define
\[\mathcal R_n=\sum_{k=0}^{K_n-1}(e^{2\pi i/K_n})^k\Chi_{A_k}\]

Then:
\begin{enumerate}
\item $\mathcal R_n$ is an eigenvector  of $U_{\mathfrak O}$ with eigenvalue 
$e^{2\pi i/K_n}$,
\item $(\mathcal R_n)^{k_n}=\mathcal R_{n-1}$,
\item $\{(\mathcal R_n)^k:0\le k< K_n, n\in\nn\}$ form a basis for $L^2(\prod_{i}\poZ_{k_i})$. 
\end{enumerate}
\medskip

For a fixed $n$, the sets $\{A_i:0\le i<K_n\}$ form a tower which will play a special role in our proofs. More generally if $\xbmt$ is an ergodic measure preserving system and $\pi:X\to \mathfrak O$ is a factor map, we set $B^i_n=\pi^{-1}A_i$. Then $\{B^i_n:0\le i<K_n\}$ is a partition of $X$ that forms a tower in the sense that $T[B^i_n]=B^{i+1}_n$ for $i<K_n-1$ and $T[B^{K_n-1}_n]=B^0_n$.

\begin{definition}\label{n-tower}
We will call the tower $\mct_n=\{B^i_n:0\le i<K_n\}$ be the \emph{$n$-tower associated with 
$\mathfrak O$}.
\end{definition}

\begin{figure}[h!]
\centering
\includegraphics[height=.35\textheight]{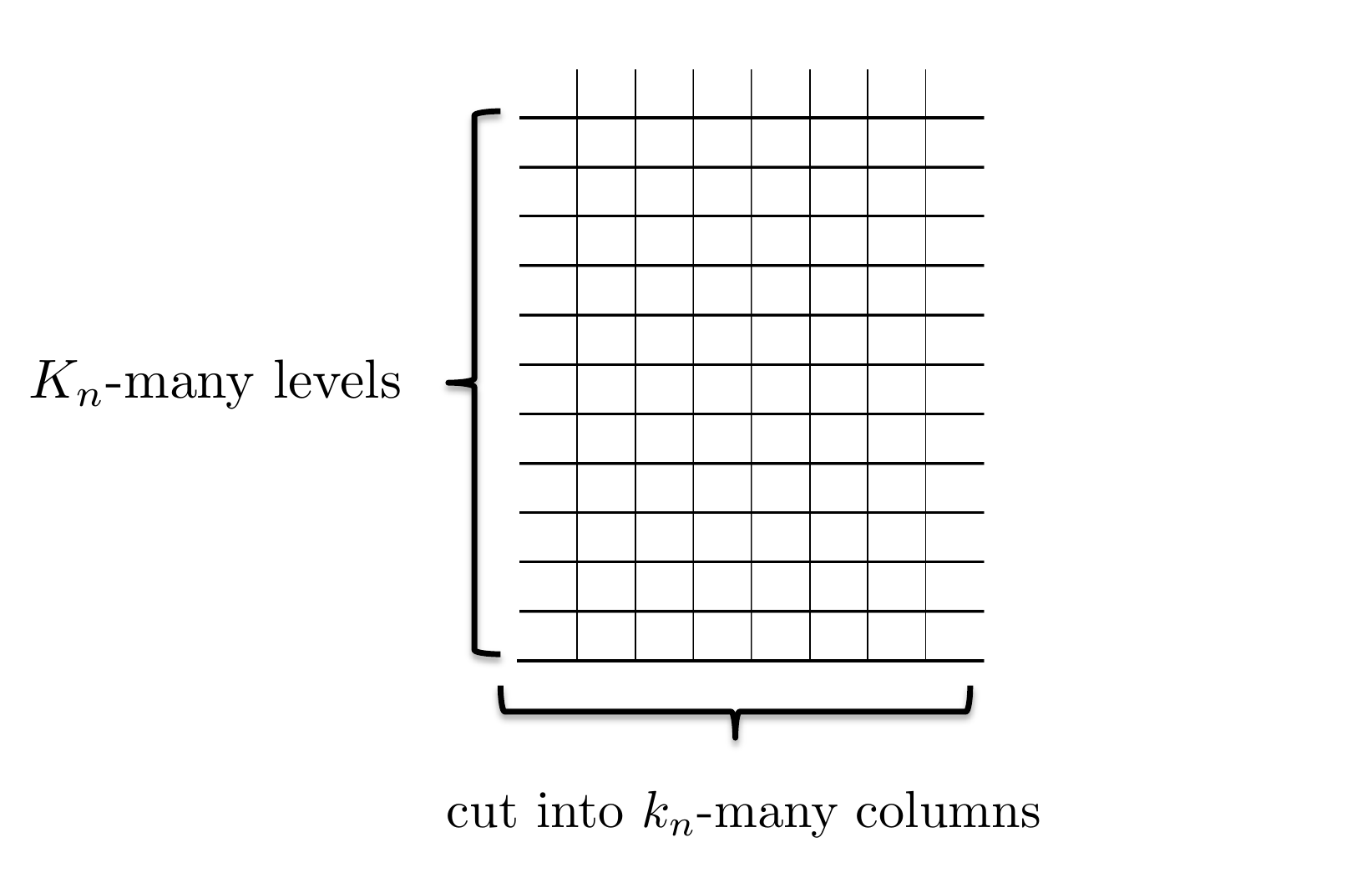}
\caption{The tower $\mct_n$.}
\label{Tncut}
\end{figure}
Figure \ref{Tncut} illustrates the $n^{th}$ tower. The horizontal lines represent the levels of the tower. The ``$n+1^{st}$-digit" of points in $\mbox{{\boldmath $O$}}$ determine $k_n$ many vertical cuts through $\mct_n$. Enumerating the levels according to their lexicographic order in $\prod_{j<{n+1}}\poZ_{k_j}$ amounts to stacking the post-cut columns of $\mct_n$:

\begin{figure}[h]
\centering
\includegraphics[height=.3\textheight]{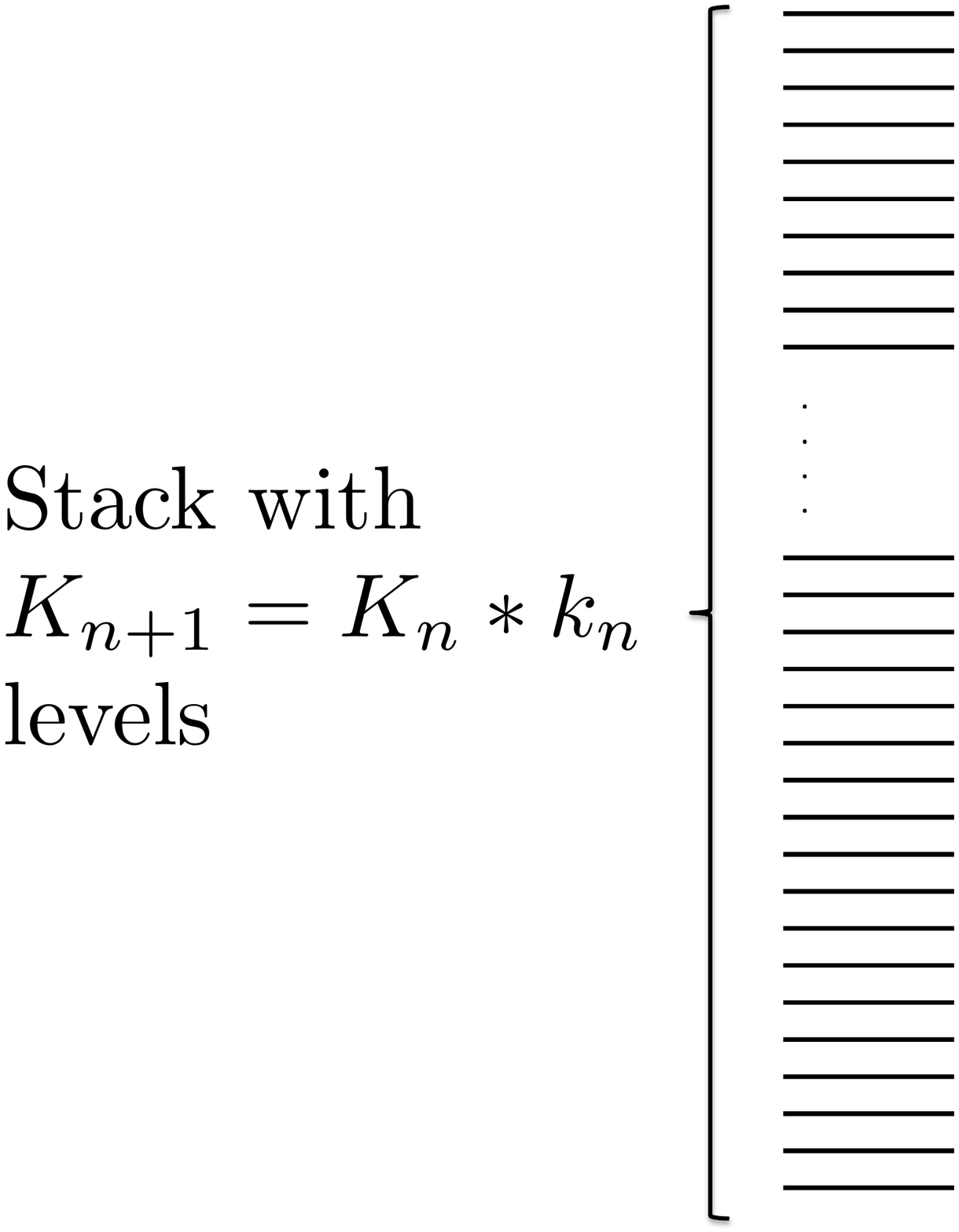}
\caption{The tower $\mct_{n+1}$. }
\label{Tn+1}
\end{figure}

\medskip

\bfni{Spectral Characterization} Here is a standard spectral characterization of transformations with an odometer factor. Suppose now that $\xbmt$ is an ergodic measure preserving system. Let $U_T:L^2(X)\to L^2(X)$ be defined by $U_T(f)=f\circ T$. Let $G$ be the group of eigenvalues of $U_T$ that have finite multiplicative order (as elements of $\mathbb C$). 

Suppose that $G$ is infinite. Then there is  a sequence of generators $\{g_n:n\in\nn\}$ of $G$ so that $o(g_{n}))|o(g_{n+1})$. The dual $\hat{G}$ of $G$ is the odometer based on $\la k_n:n\in\nn\ra$, where $k_n=o(g_n)$. We have outlined the proof of:

\begin{prop}\label{eigenvalue charact}
Let $\xbmt$ be an ergodic measure preserving system. Then $X$ has an odometer factor if and only if $U_T$ has infinitely many eigenvalues of finite  multiplicative order.
\end{prop}

Here is a useful remark. 
\begin{prop}\label{fast growing}
Let $\skn$ determine an odometer transformation $\mathfrak O$ and $K_n=\prod_{i< n}k_i$. Then for any infinite subsequence of $\la K_{n_j}:j\in\nn\ra$ of $\la K_n:n\in\nn\ra$ if we set $k'_0=K_{n_0}$ and for $j\ge 1, k'_i=K_{n_j}/K_{n_{j-1}}$, then the odometer $\mathfrak O'$ determined by 
$\la k_j':j\in\nn\ra$ is isomorphic to $\mathfrak O$. 
\end{prop}
\noindent In particular an arbitrary odometer $\mathfrak O$ has a presentation where 
$\sum 1/k_n<\infty$. 

\subsection{Invariant measures}
Let $X$ be a compact separable metric space and $T:X\to X$ be a homeomorphism.  The the space the collection of $T$-invariant probability measures on $X, \mathcal M(X,T)$, endowed with the weak topology, forms a Choquet simplex $\mathcal K$: a compact, metrizable subset of a locally convex space such that  for each $\mu\in \mathcal K$ there exists a unique measure  concentrated on the extremal points of  $\mathcal K$ which represents  $\mu$. Since the extreme points of the invariant measures are the ergodic measures, this is a statement of the Ergodic Decomposition Theorem. 


 \section{Odometer Based Symbolic Systems}\label{OBSS}
Here is the general definition of a construction sequence and its limit. We will be working with a special case, the odometer construction sequences.
 \begin{definition}\label{constseq}
 A \emph{construction sequence} in a finite alphabet $\Sigma$ is  
 a sequence of collections of words $\la\mcw_n:n\in\nn\ra$ with the properties that:
\begin{enumerate}
\item $\mcw_0=\Sigma$,
\item all of the words in each $\mcw_n$ have the same length $q_n$ and are uniquely readable,
\item each $w\in\mcw_{n}$ occurs  at least once as a subword of every $w'\in \mcw_{n+1}$,\label{minimal?}

\item \label{not too much space} there is a summable sequence $\la \epsilon_n:n\in\nn\ra$ of positive numbers such that for each $n$, every word $w\in \mcw_{n+1}$ can be uniquely parsed into segments 
\begin{equation}u_0w_0u_1w_1\dots  w_lu_{l+1}\label{words and spacers}
\end{equation}
 such that each $w_i\in \mcw_n$, $u_i\in \Sigma^{<\nn}$ and for this parsing
\begin{equation}\label{small boundary numeric} {\sum_i|u_i|\over q_{n+1}}<\epsilon_{n+1}.
\end{equation}
\end{enumerate}
 \end{definition}
\noindent We call the elements of $\mcw_n$ ``$n$-words," and let $s_n=|\mcw_n|$.

\begin{definition} Let $\bk$  be the collection of 
 $x\in \Sigma^\poZ$ such that every finite contiguous subword of $x$ occurs inside some $w\in \mcw_n$. Suppose $x\in \bk$ is such that   $a_n\le 0< b_n$  and $x\rest[a_n, b_n)\in \mcw_n$. Then $w=x\rest[a_n,b_n)$ is the \emph{principal $n$-subword} of $s$. We set $r_n(s)=|a_n|$, which is the position of $s(0)$ in $w$. 
 \end{definition}
\noindent Then $\bk$ is a 
 closed shift-invariant subset of $\sz$ that is compact if $\Sigma$ is finite. Clause 3.) of the definition guarantees that $\bk$ is indecomposable as a topological system.
 
 Not every symbolic shift in a finite alphabet can be built as a limit of a construction sequence, however this method directly codes cut-and-stack constructions of transformations on probability spaces.

\medskip

Odometer construction sequences are those that use no spacers $u_i$:
\begin{definition}\label{distinguishing odos}Let $\la k_n:n\in\nn\ra$ be a coefficient sequence,
A construction sequence $\la \mcw_n:n\in\nn\ra$  is  \emph{odometer based} if and only if  $\mcw_{n+1}\subseteq \mcw_n^{k_n}$. 
 A symbolic system $\bk$ is \emph{odometer based} if it has a construction sequence that is odometer based.  For an odometer based construction sequence we let $K_n=\prod_{m<n}k_m$.\footnote{$K_n$ will be equal to the $q_n$ in definition \ref{constseq}.}
 \end{definition} 
 \noindent For odometer based construction sequences,  strengthening clause 3.) in definition \ref{constseq} of \emph{construction sequence} to require that each $w\in \mcw_n$ occurs at least twice in every $w'\in \mcw_{n+1}$ has the consequence that $\bk$ is a minimal system.
\smallskip

In this section we prove 

\begin{theorem}\label{first presentation}
 Let $\xbmt$ be a measure preserving system with finite entropy. Then $X$ has an odometer factor if and only if $X$ is isomorphic to a topologically minimal odometer based symbolic system.
\end{theorem}

 If $\bk$ is an odometer based system with construction sequence $\la \mcw_n:n\in\nn\ra$, then for all $s\in \bk$ there are $a_n\le 0\le b_n$ such that 
 $s\rest [a_n, b_n)\in \mcw_n$. In particular for every $n$, every $s\in \bk$ has a principal $n$-subword.

 The name \emph{odometer based system} is motivated by the following proposition:
 
 \begin{proposition}\label{canfact} Suppose that $\la \mcw_n:n\in\nn\ra$ is an odometer based construction sequence for a symbolic system $\bk$. Let $K_n$ be the length of the words in $\mcw_n$, $k_0=K_1$ and for $n>0$,  $k_n=K_{n+1}/K_{n}$. Then the odometer $\mathfrak O$ determined by $\la k_n:n\in\nn\ra$ is canonically a factor of $\bk$.
 \end{proposition}
 \pf Let $s\in \bk$. By the unique readability, for each $n$, $s$ can be uniquely parsed into a bi-infinite sequence of $n$ words. For each $n$, there is an $c_n$ such the principal $n$-block of $s$ is the   $c_n^{th}$ $n$-word in the principal $n+1$-block of $s$. 
 
 Define a map $\phi:\bk\to \prod_{n} \poZ/k_n\poZ$ by setting $\phi(s)=\la c_n:n\in\nn\ra$. It is easy to check that $\phi(sh(s))=\mco(\phi(s))$.
 \qed
 
 One way of defining elements of $\bk$ is illustrated in the following Lemma. 
 \begin{lemma}\label{defining s}
 Let $\la r_n:n\ge k\ra$ be a sequence of natural numbers and $\la w_n:n\ge k\ra$ be a sequence of words with $w_n\in \mcw_n$.  Suppose that for each $n$, the $r_n^{th}$ letter in $w_{n+1}$ is inside an occurrence of $w_n$ in $w_{n+1}$. Then there is a unique $s\in \bk$ such that for $n\ge k$,  $r_n(s)=r_n$ and the  principal $n$-subword of $s$ is $w_n$.
 \end{lemma}
 \bigskip

 \subsection{Odometers are not topological subshifts}
Theorem \ref{first presentation} says that all ergodic measure preserving transformations with a non-trivial odometer factor are measure theoretically isomorphic to an odometer based symbolic system. In contrast, it is well known that as \emph{topological} dynamical systems, odometers are not homeomorphic to symbolic shifts. For background we give a very brief proof of this fact.

 \begin{definition}  Let $(X,d)$ be a metric space. A map $T:X\to X$ is \emph{expansive} if there is an $\epsilon>0$ such that for all $x\ne y$ in $X$ there is an $n, d(T^nx,T^ny)\ge\epsilon$.
\end{definition}

The following is easy to verify:
\begin{proposition}\label{subshiftsexpansive}
Let $(X,d)$ be a compact metric space and $T:X\to X$.
\begin{enumerate}
\item If $T$ is an isometry then then $T$ is not expansive unless $X$ is finite.
\item If $X\subseteq \Sigma^\poZ$ is a compact subshift, and $T$ is the shift map, then $T$ is expansive. 
\end{enumerate}
\end{proposition}
\pf The first proposition is trivial. To see the second, note that we can assume $\Sigma$ is finite.  Let $c$ be the minimum distance between cylinder sets $\la i\ra$ and $\la j\ra$ based at $0$. Then if $x\ne y$, we can find an $n, x(n)\ne y(n)$. It follows that $d(T^nx,T^ny)\ge c$.\qed

 In view of Proposition \ref{subshiftsexpansive}, to see that an odometer cannot be presented as a topological subshift it suffices to show that, viewed as metric systems,   odometer transformations are isometries. Let $O=\prod_0^\infty \poZ/k_n\poZ$ be an odometer and $T$ be the  odometer map 
 $\mco$.
 
 For $x, y\in O$, define $\Delta(x,y)$ to be the least $n$ such that $x(n)\ne y(n)$ and $d(x,y)={1\over 2^{\Delta(x,y)}}$. Then $d$ is a complete metric yielding the product topology on $O$ and is invariant under $\mco$. Thus, by Proposition \ref{subshiftsexpansive}, if follows the  that odometer is not isometric to a subshift of $\Sigma^\poZ$ for any finite $\Sigma$.
  
  \subsection{Presenting the Odometer}\label{just the odo}
To illustrate one of the main ideas in the proof we give a presentation of an arbitrary odometer as an odometer based system.

\begin{example}\label{odosolo}
If $\mathfrak O$ is an odometer determined by $\la k_n:n\in\nn\ra$ with $k_n\ge 2$, then there is an odometer based construction sequence $\la \mcw_n:n\in\nn\ra$ such that the associated symbolic system $\bk$ is uniquely ergodic and measure theoretically conjugate to $\mathfrak O$. 
\end{example}  
 \pf By Proposition \ref{fast growing}, we can assume that $\sum 1/k_n<\infty$.  We define an odometer based construction sequence $\la \mcw_n:n\in\nn\ra$ such that each $\mcw_n=\{a_n, b_n\}$ has exactly two words in it.
 \begin{itemize}
 \item Let $\Sigma=\{a,b\}$ and $\mcw_0=\Sigma$. 
  \item Suppose that we are given $\mcw_n=\{a_n, b_n\}$. Let $\mcw_{n+1}=\{a_{n+1}, b_{n+1}\}$ with $a_{n+1}, b_{n+1}\in \mcw_n^{k_n}$ where:
  \begin{eqnarray*}
  a_{n+1}&=&a_na_na_nb_nb_nb_na_nb_na_nb_n\dots x\\
  b_{n+1}&=&b_nb_nb_na_na_na_na_nb_na_nb_n\dots x
  \end{eqnarray*}
where $x$ is either $a_n$ or $b_n$, depending on whether $k_n$ is  even or odd.
\end{itemize}
It is easy to verify inductively that the the $a_n$'s and $b_n$'s are uniquely readable (look for patterns of the form $a_na_na_n$ and $b_nb_nb_n$) and that $\la\mcw_n:n\in\nn\ra$ is uniform. Let $\bk$ be the associated symbolic system. Then $\bk$ is uniquely ergodic, with an invariant measure $\mu$. 

Let $\phi:\bk\to \mco$ be the canonical map from Proposition \ref{canfact}. To establish the claim in the Example \ref{odosolo} it suffices to show that there is a set of measure one for the odometer on which $\phi$ is invertible. 

Let $G=\{x\in \mco:$ for all large enough $n, x(n)\ge 10\}$. Since $\sum 1/k_n<\infty$, the Borel-Cantelli Lemma implies that  $G$ has measure one for 
$\mathfrak O$. 

We define $\psi:G\to \bk$ so that $\phi\circ\psi=id$. By Lemma \ref{defining s}, we can determine $\psi(x)$ by defining a suitable sequence 
$\la r_n:n\ge k\ra$ and $\la w_n:n\ge k\ra$. 

Let $x\in G$ and suppose that for all 
$n\ge k, x(n)\ge 10$. Fix $n\ge k$.
\[r_n=x(0)+x(1)k_0+x(2)k_1+\cdots+x(n)k_n.\]
Since $x(n)\ge 10$, either for all $n+1$-words $w\in \mcw_{n+1}$,  the $x(n)^{th}$ $n$-subword in 
$w$ is $a_n$ or for all  $n+1$-words $w\in \mcw_{n+1}$, the $x(n)^{th}$ $n$-subword in $w$ is 
$b_n$. Let $w_n$ be either $a_n$ or $b_n$ accordingly.

Let $\psi(x)$ be the element $s$ of $\bk$ determined by $\la r_n:n\ge k\ra$ and $\la w_n:n\ge k\ra$. Then $\psi(x)$ is well-defined and $\phi\circ\psi(x)=id$.  If $\nu$ is the  measure on $\mbox{{\boldmath $O$}}$ giving the odometer system,  then $\psi$ induces a shift-invariant measure $\nu^*=\psi^*\nu$ on 
$\bk$.  Since $\bk$ is uniquely ergodic, $\nu^*=\mu$ and $\psi=\phi^{-1}$.\qed
\medskip
We note that the set $G$ in the proof is a Borel set and $\psi$ is continuous.
\subsection{The plan}\label{panama}

In this section we explain the idea of the proof of Theorem \ref{first presentation}, the details will follow in the next section.
 To show that a given transformation with an odometer factor is isomorphic to a symbolic system built from an odometer based construction sequence we build a generating partition so that the  names of points on the bases of the $n$-towers in Definition \ref{n-tower} form an odometer based construction sequence.

Let $\xbmt$ be an ergodic measure preserving system with an odometer factor $\mathfrak O$.  By Example \ref{odosolo}, $\mathfrak O$ is isomorphic to an odometer based system in the alphabet $\Sigma=\{a,b\}$. Call the resulting construction sequence $\la \mcw_n^\mco:n\in\nn\ra$. If $\bk$ is the symbolic system associated with this construction sequence  we have:
\[
\begin{diagram}
\node{X}\arrow[1]{e,r}{\pi}\node{\mco}\arrow{e,r}{\phi}\node{\bk}
\end{diagram}
\]
Let  $\mcq=\{Q_0, Q_1\}$ be the partition of $X$ corresponding to the basic open intervals $\la a\ra, \la b\ra$ in $\bk$ (so $Q_i=(\phi\circ\pi)\inv\la i\ra$). Then 
$\mcq$ generates the factor 
$\mathfrak O$.

Suppose that $C\subseteq X$ is a set of positive measure. Let $T_C:C\to C$ be the induced map: $T_C(c)=d$ if and only if for the least $k>0, T^k(c)\in C$ one has $T^k(c)=d$. Suppose that $\mcp_0=\{P_1, P_2, \dots P_a\}$ is a generator for $T_C$, where $a\in \nn$, $D=X\setminus C$ and $\mcp=\mcp_0\cup\{D\}$.
Then for $x\in X$, the $\mcp$-name of $x$ uniquely determines $x$, and thus $\mcp$ is a generator for $X$.

\begin{figure}[h]
\includegraphics[height=.3\textheight]{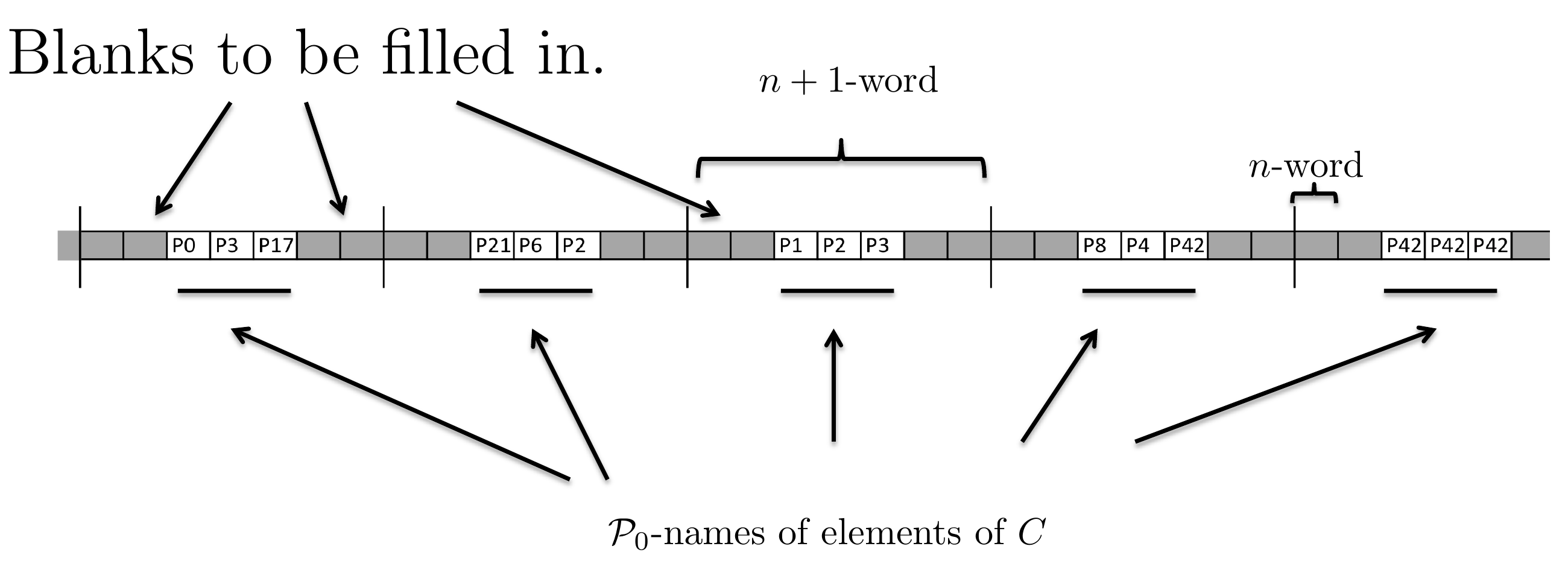}
\caption{The $\mcp_0$-name of $x$ punctuated by the odometer. }
\label{x-orbit}
\end{figure}

 For a typical $x$, the combined $\mcp_0, \mcq$-name of $x$ can be visualized as in figure \ref{x-orbit}.
The elements of $\mcq$ parse the $x$-orbit into $n$-words which determine the duration an orbit stays in $D$, while the elements of $\mcp_0$ determine the orbit of $x$ inside $C$. Since $\mcp_0$ and $\mcq$ determine $x$, in building an odometer based symbolic representation of $\xbmt$, one has complete freedom to fill in symbols in the parts of the $x$-orbit that  lie in $D$.  This allows our word construction to satisfy the definition of \emph{odometer based}.

In terms of partitions, this can be restated as saying that we can modify the atoms  of the partition  $\mcp_0$ by adding elements of $D$ in any arbitrary way, as long as the restriction of each atom of $\mcp_0$ to $C$ remains the same. If $\mcp_0'=\{P_1', P_2', \dots P_a'\}$ is the modification of $\mcp_0$, then any partition refining  $\mcp_0'$ and $\mcq$ still forms a generator for $T$. Hence, as in Remark \ref{joining factors}, the symbolic system consisting of pairs $(s_{\mcp_0'}, s_\mcq)$ of $\mcp_0'$ and $\mcq$-names is isomorphic to $\xbmt$.

Of course, Figure \ref{x-orbit} is an over-simplification of the possibilities for the orbit: it assumes that the set $C$ fits coherently with the odometer factor. In other words,  $C$ must be chosen to be measurable with respect to the sub-$\sigma$-algebra of $\mcb$ generated by the odometer factor.

\subsection{The proof}
Suppose that $\xbmt$ has entropy less than ${1\over 2}\log a$. By Proposition \ref{fast growing},  we can assume that $K_1=k_0>10$, $K_n=\prod_{i<n}k_i$, and $k_{n}>4a^{K_n}10^{n+1}$.

Let $B_0, B_1, \dots$ be the bases of the $n$-towers in $X$ associated with 
$\mathfrak O$; in the notation of Definition \ref{n-tower}, $B_n=B^0_n$. Let $d_n=4K_{n-1}a^{K_{n-1}}$ and define
\begin{eqnarray*}
D_n&=&\bigcup_{0\le i\le d_n}B_n^i\\
&\mbox{and}&\\
D&=&\bigcup_1^\infty D_n
\end{eqnarray*}
Thus $D_n$ consists of the first $d_n$ levels of the $n$-tower. Since all of the levels of the tower have the same measure the measure of $D_n$ is
\begin{eqnarray*}
{d_n\over K_n}&=&{4K_{n-1}a^{K_{n-1}}\over K_n}\\
&=&{4K_{n-1}a^{K_{n-1}}\over K_{n-1}k_{n-1}}\\
&<&{4K_{n-1}a^{K_{n-1}}\over K_{n-1}4a^{K_{n-1}}10^{n+1}}\\
&=&10^{-(n+1)}.
\end{eqnarray*}

Set $C=X\setminus D$. Clearly $C$ is measurable with respect to the odometer factor, since it is a union of levels of the odometer towers. Moverover, 
 $\mu(C)>3/4$, and hence the entropy of $T_C$ is less than $(2/3)\log a$. By Krieger's Theorem \cite{Krieger} there is a generating partition $\mcp_0=\{P_1, P_2, \dots P_a\}$ for $T_C$, that has $a$ elements. We can assume without loss of generality that $a\ge 2$.

\begin{figure}[h]
\centering
\includegraphics[height=.6\textheight]{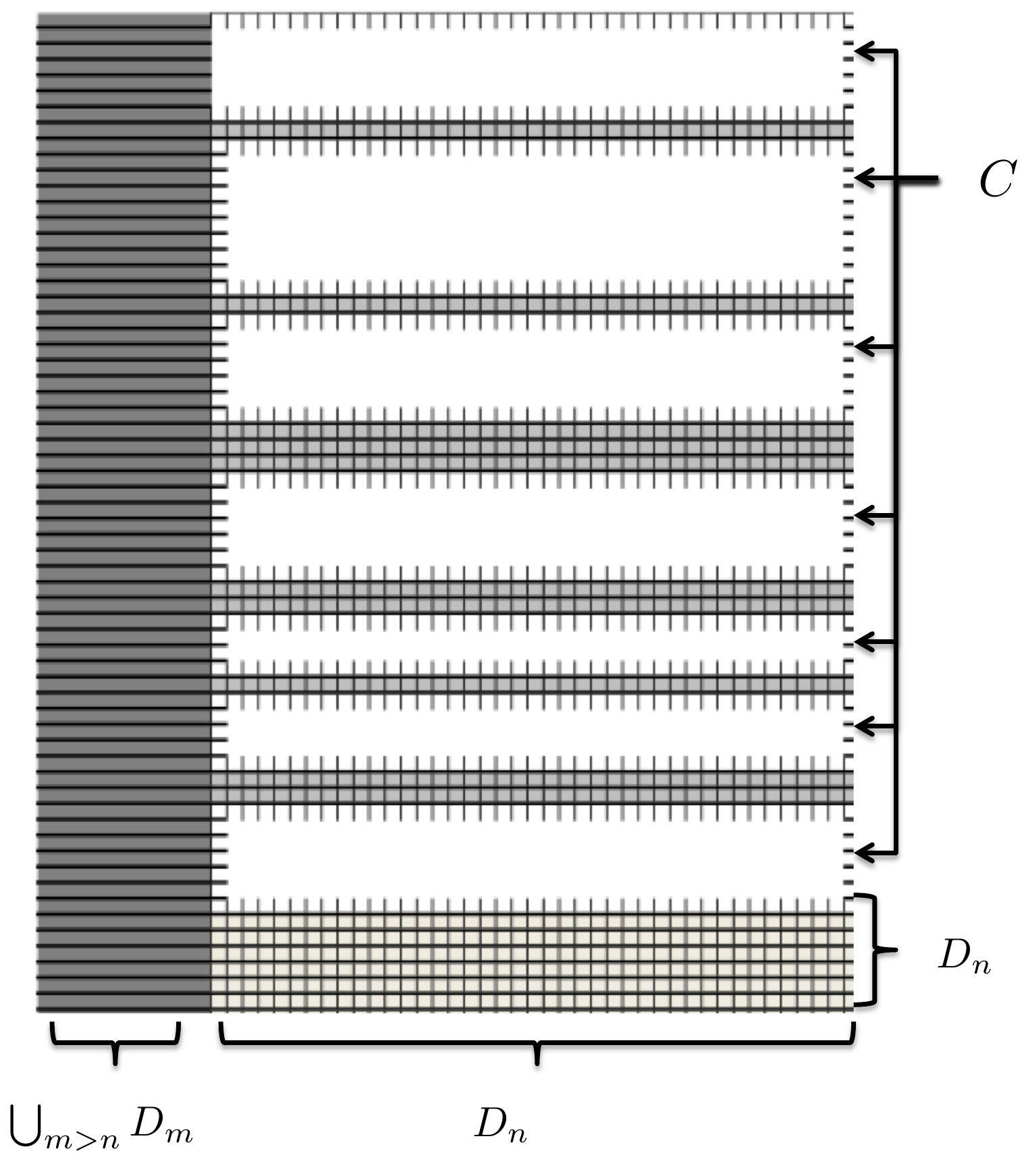}
\caption{The $n^{th}$ stage of the construction. The shaded horizontal bands are elements of $D_m$ for $m<n$.}
\label{core}
\end{figure}

Figure \ref{core} is a graphical representation of $\mct_n$ showing:
\begin{enumerate} 
\item $C$ as whitespace, 
\item $D_n$ lightly shaded as an initial segment of the levels of $\mct_n$
\item The sets $D_m$ for $m<n$ are initial segments of earlier $\mct_m$ and hence get stacked as bands across $\mct_n$. They are given an intermediate shading in figure \ref{core}.
\item 
Because each  $D_m$ is an  initial segment of $\mct_m$, at the previous stage the points in $D_{m-1}$ have to be in the leftmost columns of $\mct_{m-1}$. Moreover for $m<m', K_m$ divides $d_{m'}$. Thus $D_{m'}$ is made up  of whole columns of $\mct_m$.  Consequently $\bigcup_{m>n}D_m$ forms a contiguous rectangle on the left side of $\mct_n$. This region is indicated by the darkest shading.
\end{enumerate}

We construct $\mcp_0'$ in the manner described in Section \ref{panama}: we add points from $D$ to each $P_j$ to get a final partition $\mcp_0'=\{P_1', P_2', \dots P_a'\}$. The corresponding construction sequence will use  the alphabet $\Sigma=\{P_1', P_2' \dots P_a'\}\times \mcq$. $\mcw_0=\Sigma$ and 
$\mcw_n$ will consist of the $\Sigma$-names of points  that occur in the base of $\mct_n$.  Thus the construction is completely determined by the manner we add points  to the $P_i$.

The words must satisfy Definition \ref{constseq}. Clause 1 is automatic.
 Clause 2 holds because all words have length equal to the the height of 
 $\mct_n$. Unique readability is immediate since the odometer based presentation of $\mathfrak O$ uses uniquely readable words in the language 
 $\mcq$. Clause 4 is vacuous since we have no spacers $u_i$ occurring anywhere in the words: elements of $\mcw_{n+1}$ are simply concatenations of words from $\mcw_n$. 
 
 The system is minimal if each word in $\mcw_n$ occurs at least twice 
 in each word in $\mcw_{n+1}$, a property which is stronger than clause 3. We satisfy this by ``painting" the words from $\mcw_n$ onto $D_{n+1}$.

 Let $P'_i(n)$ be the collection of points in $P_i'$ at stage $n$, and $P'_i(0)=P_i$ Inductively we will assume that at stage $n$:
 \begin{enumerate}
 \item    
 $\bigcup_{n<m}D_m\cap P'_i(n)=\emptyset$ for all $i$
 and
 \item $(\mct_n\setminus \bigcup_{n<m}D_n)\subseteq \bigcup_iP'_i(n)$.
 \end{enumerate}

For $n=1$, we consider $D_1\setminus\bigcup_{m>1}D_m$. At stage 1 the minimality requirement says  that each pair $(P_i'(0),j)$ for $1\le i\le a$ and $j\in \{0,1\}$ occurs at least {twice}.  Each of $0,1$ occur equally often in the $\mcq$-names of the first $d_1$ letters of each $\mcq$-name and $d_1=4a$. Hence it is possible to assign the levels in $D_1\setminus \bigcup_{m>1}D_m$ to $\{P_1'(1), \dots P_a'(1)\}$ in such a way that each  $(P_i'(0),j)$ occurs at least twice.

To pass from $n$ to $n+1$ in the construction, we know inductively that no elements of $D_{n+1}$ have been assigned to any $P_i'$ at earlier stages. 
Moreover $\mcw_n$ consists of the $\Sigma$-names of the words in 
$B_0\setminus \bigcup_{m>n}D_m$, where $B_0$ is the base of $\mct_n$. There are at most $2a^{K_n}$ such words in the language $\Sigma$. Each such word has length $K_n$.

Since $d_{n+1}=4K_na^{K_n}$ there are ample levels in $D_{n+1}$ that 
each level can be added to some $P_i'(n+1)$ in a manner that each word in 
$\mcw_n$ occurs at least twice as a $\Sigma$-name of an element the first $d_{n+1}$ levels of $\mct_{n+1}$.
\qed
\begin{remark}\label{shifting sands}
The construction in the proof of Theorem \ref{first presentation} was used a particular presentation of $\mathfrak O$ as an odometer based system in a language $\mcq=\{a,b\}$ to build a language $\Sigma=\{P_1', P_2' \dots P_a'\}\times \mcq$. If we were given another odometer based presentation $\la \mcw^{\mathfrak O}_n:n\in\nn\ra$ of $\mathfrak O$ in a different finite language with letters $\{a_1, \dots, a_k\}$ we could take $\Sigma=\{P_1', P_2' \dots P_a'\}\times \{a_1, \dots,  a_k\}$ and repeat the same construction over this presentation.  We will call this the odometer based presentation of $X$ \emph{built over $\la \mcw^{\mathfrak O}_n:n\in\nn\ra$}.
\end{remark}

\section{Toeplitz Systems}\label{Toeplitz Systems}
 
 In this section we use a result of Downarowicz (\cite{downer}) to show that every compact metrizable Choquet simplex is affinely homeomorphic to the simplex of invariant measures of an odometer based system.  Williams showed that the orbit closure of every Toeplitz sequence in a finite language $\Sigma$ is a minimal symbolic shift $\bl$ with a continuous map to an odometer factor $\mathfrak O$.  If $\pi: \bl \to \mathfrak O$ is this factor map,   it would be tempting to argue that the words occurring on $\pi$-pullbacks of the levels of the $n$-towers form an odometer based construction sequence. However we don't know this in general; in particular we don't know that the words constructed this way are uniquely readable.
 
To make the words uniquely readable we need to make the map $\pi$ ``extremely Lipschitz." To do this we introduce the  \emph{ad hoc} notion of an augmented symbolic system. 
 
 \begin{definition}
 Let $X$ and $Y$ be minimal symbolic topological shifts in alphabets $\Sigma, \Gamma$.  An \emph{augmentation} of $X$ by $Y$ is a shift-invariant Borel set $A\subseteq X\times Y$ such that if 
  $L=\{x:$ there is exactly one $y, (x,y)\in A\}$, then for all shift-invariant $\mu$ on $X$, $\mu(L)=1$. 

 \end{definition}
 \noindent We write $X|Y$ for an augmentation of $X$ by $Y$. 
 \smallskip

 Consequently:
 \begin{prop} \label{why bother?}
 Suppose that $Y$ is uniquely ergodic and $X|Y$ is an augmentation. Then there is a canonical affine homeomorphism of $\mathcal M(X,sh)$ with $\mathcal M(X|Y, sh)$.
 \end{prop}
 \pf  If $\mu$ is a measure on $X$ then $\mu$ determines a measure on $L$ and hence on $X|Y$.  Conversely if $\nu$ is a measure on $X|Y$, let $\mu=\nu^X$.  Then, since $\mu(L)=1$, 
 $\nu(\{(x,y)\in A:x\in L\})=1$ and for $B\subseteq A$, $\nu(B)=1$ if and only if $\mu(\pi_X(B))=1$. 
Thus there is  a bijection between $\mathcal M(X,sh)$ and $\mathcal M(X|Y, sh)$ that is easily seen to be an affine homeomorphism.\qed

 To prove Proposition \ref{after downer}, we use:
\medskip
 
 \bfni{Theorem}(Downarowicz, \cite{downer}, Theorem 5) For every compact metric Choquet simplex $K$ there is a dyadic Toeplitz flow whose set of invariant measures is affinely homeomorphic to $K$.
 \medskip

 \begin{prop}\label{after downer}Let $\bl$ be the orbit closure of a Toeplitz sequence $x$, $\mathfrak O$ be its maximal odometer factor based on  a sufficiently fast growing sequence $\la k_n\ra$ and $\bk$ be the  odometer based presentation of $\mathfrak O$ defined in example \ref{odosolo}.  Then there is an odometer based system $\bl^*\subseteq \bl\times \bk$ and a set $A\subseteq \bl^*$ that is an augmentation of $\bl$ by $\bk$ and has measure one for every invariant measure on $\bl^*$.
 \end{prop}

 Thus, as an immediate consequence of Downarowicz' theorem and Propositions \ref{why bother?}, \ref{after downer}:
 
 \begin{corollary}\label{much ado}
For every compact metrizable Choquet simplex there is an odometer based symbolic shift $\bl^*$ whose set of invariant measures is affinely homeomorphic to $K$. 
 \end{corollary}
 
\pf(Proposition \ref{after downer})
We use the language of Williams \cite{williams}.
Let $x$ be a Toeplitz sequence in a finite language $\Sigma$.  Let $\bl$ be the orbit closure of $x$ under the shift map and $\mathfrak O$ be the associated 
 odometer system. 
 
As in \cite{williams} we can choose a sequence $\la K_n:n\in \nn\ra$ of essential periods for $x$.  By choosing the $K_n$'s to grow fast enough we can assume that 
    \begin{enumerate}
    \item[a.)] $K_n|K_{n+1}$
    \item[b.)] $\bigcup_{n}$Per$_{K_n}(x)=\poZ$.
    \end{enumerate}
    
 Choosing a further subsequence we can also assume
  that
    \begin{enumerate}
    \item[c.)]  if $k\equiv0$(mod $K_n$) then there is an $i\equiv 0$(mod $K_n$) with $i<K_{n+1}$ and $x\rest [k, k+K_n)=x\rest[i, i+K_n)$.
    \end{enumerate}
 Given $n_0$, for large enough $n, x\rest[0, K_{n_0})$ is a subset of the $K_n$-skeleton of $x$. Since the $K_n$-skeleton is $K_n$-periodic, every subword of the $K_n$-skeleton is repeated $K_{n+1}/K_n$ times in $x\rest[0,K_{n+1})$. Thus by again thinning the $K_n$'s we can assume that:
     \begin{enumerate}
     \item[d.)] for each $n$ and $i\equiv0$(mod $K_n$) and each word $w\in \Sigma^{K_n}$ occurring as $x\rest[i, i+K_n)$, $w$ occurs at least twice in $x\rest[0, K_{n+1})$.
     \end{enumerate}
     
 Let $\mathfrak O$  be the odometer with coefficient sequence $\la k_n:n\in\nn\ra$, where $k_n=K_{n+1}/K_n$.
Let $\la \mcw_n:n\in\nn\ra$ be the odometer based construction sequence in the presentation of $\mathfrak O$ given in Example \ref{odosolo}. Let $w^0_n, w^1_n$ be the two words in $\mcw_n$. We define an odometer based construction sequence by setting $\mcv_n$ to be the collection of words $v$ in the alphabet $\Sigma\times \{a,b\}$ of the form 
\[(x\rest(i, i+K_n), w^j_n)\]
where $i<K_{n+1}, i\equiv0$(mod $K_n$) and $j\in \{0,1\}$.

To see that this is an odometer based construction sequence we check definitions \ref{constseq} and \ref{distinguishing odos}.

Unique readability of the words $v\in \mcq_n$ follows immediately from the fact that the $w^j_n$ are. The fact that each $w\in \mcw_n$ occurs at least once as a subword of each $w'\in \mcw_{n+1}$ follows immediately from item d.) of the properties of the essential periods of $x$. From item c.) and  structure of the word construction each word in $\mcv_{n+1}$ is a concatenation of words in $\mcv_n$.

By \cite{williams}, there is a continuous factor map 
	\[\pi:\bl\to \mathfrak O.\] 
From Example \ref{odosolo} we see that there is an invariant Borel set $G\subset \mathfrak O$ of measure one and a one-to-one, continuous map 
$\psi:G\to \mathbb K$. Let  $\bl^*$ be the limit of this construction sequence and
\[A=\{(y,\psi\circ\pi(y)):\pi(y)\in G\}\subseteq \bl^*.\]
Let $\mu$ be an invariant measure on $\bl$. Then $\mu(\pi^{-1}(G))=1$, and for $y\in \pi^{-1}(G)$ there is a unique $z, (y,z)\in A$.

Let $\rho$ be an invariant measure on $\bl^*$. Let $\rho^\bl$ be the $\bl$ marginal. Then $\rho^\bl(\pi^{-1}(G))=1$. If $y\in\pi^{-1}(G)$ and $(y,z)\in \bl^*$, then $z=\psi\circ\pi(y)$.  Hence {$\mu(A)=1$.}\qed

 \noindent The next example is an odometer based system that is far from being a Toeplitz system.
\begin{example}\label{aperiodic odometer based}
There is an odometer based system $\bk$ such that no $x\in \bk$ has any periodic locations: for all $x\in \bk, p\in \nn,  \mbox{Per}_p(x)=\emptyset$. In particular no 
$x\in \bk$ is a Toeplitz sequence.
\end{example} 
 \pf   Let $\Sigma=\{0,1\}$. For $w\in \Sigma^{<\nn}$ define $\bar{w}$ to the result of substituting $0$'s for the $1$'s in $w$ and vice versa.
 
 Define an odometer based construction sequence $\la\mcw_n:n\in\nn\ra$ by induction. Let $\mcw_0=\{0,1\}$. At stage $n+1$ we will assume that each $\mcw_n$ is of the form 
 $\{w, \bar{w}\}$ where $w$ has length $K_n$. 
 Let $v=w^{K_n}\bar{w}^{K_n}$ and $\mcw_{n+1}=\{v,\bar{v}\}$. We note that 
 $\bar{w}^{K_n}w^{K_n}=\bar{v}$ so this description is unambiguous. 
 \medskip
 
 \bfni{Claim} Let $\bk$ be the symbolic system associated with $\la \mcw_n:n\in\nn\ra$. Then  for all $x\in \bk$, $k\in \poZ$, $p\in\nn$ there is a $ b\in \poZ$ 
 such that $x(k)\ne x(k+bp)$. 
 
 \pf Fix $x$, $k$ and $p$. Let $n$ be so large that $k$ and $k+p$ are in the principal $n$-block of $x$. Let $w$ be the principal $n$-subword of $x$ and assume first that the principal $n+1$-subword of $x$ is of the form $v=w^{K_n}\bar{w}^{K_n}$. Since $p<K_n$ there is an $a>0$ such that $k+apK_n\in [K_n^2,2K_n^2)$. Let $b=aK_n$. Then $[k]_{K_n}=[k+bp]_{K_n}$ and the $(k+bp)^{th}$ position of $v$ is in $\bar{w}$. It follows that $x(k)\ne x(k+bp)$. 
 
 The case where the principal $n$-subword of $x$ is in the second half of the principal $n+1$ subword is the same,  except that $a<0$.
 \qed

 We note we have proved something much stronger than claimed in the statement of Example \ref{aperiodic odometer based}, namely in the notation of \cite{williams}, for all $x\in \bk, \sigma\in \Sigma$ we have 
 Per$_p(x,\sigma)=\emptyset$.

\section{The small word property and rates of descent}\label{SWP}

The applications of the representation theorem and Proposition \ref{after downer} require that for all invariant measures on the limiting system $\bk$, the basic open intervals determined by words in 
$\mcw_{n+1}$ have measure much smaller than the measures of basic open intervals determined by words in $\mcw_n$. We show how to arrange this for odometer based systems by taking subsequences.

We define the frequency of occurrences of $w$ in $w'$, to be 
\[Freq(w,w')=\frac{ \mbox{number of occurrences of $w$ in $w'$}}{K_m/K_n}.\]
For $n<m$, 
clause 3 of the definition of a construction sequence (Definition \ref{constseq}) implies that the frequency of each word $w\in \mcw_n$ inside each $w'\in \mcw_{m}$ is at least $1/k_n$.

\begin{remark}\label{hurry up} Let $w\in \mcw_k$.
If for all $w'\in \mcw_{k+1}$, $\eta_0<Freq(w,w')<\eta_1$,   then for $k+l>k, w'\in \mcw_{k+l}$ we have $\eta_0<Freq(w,w')<\eta_1$. 
\end{remark}

\begin{definition}\label{small words} Let $\la\mcw_n:n\in\nn\ra$ be an odometer based construction sequence. Let $f_n=\sup\{Freq(w,w'):w\in \mcw_n, w'\in\mcw_{n+1}\}$ be the supremum of the frequencies of the $n$-words in $n+1$-words. The sequence $\la \mcw_n:n\in\nn\ra$ has the \emph{small word property} with respect to a sequence 
  $\la \delta_n:n\in \nn\ra$ if and only if 
  for all $n$ $f_n<\delta_n$.
  \end{definition}

The next lemma follows immediately from the Ergodic Theorem:
\begin{lemma}\label{small but not zero} Let $\la \mcw_n:n\in\nn\ra$ be an odometer based construction sequence for the system $\bk$, and $\rho$ be a shift-invariant measure  on $\bk$. Then for all words $w\in \mcw_n$:
\begin{equation*}{1\over K_{n+1}}\le \rho(\la w\ra)\le {f_n\over K_{n}}.
\end{equation*}
\end{lemma}

Thus if $\la \mcw_n:n\in\nn\ra$ has the small word property with respect to $\la \delta_n:n\in\nn\ra$ with $\delta_n<1$ then for all $w\in \mcw_n, w'\in \mcw_{n+1}$ and all invariant measures $\rho$:
\begin{equation}
\label{the point of swp}
\rho(\la w'\ra)<{\delta_{n+1}\over K_{n+1}}<\rho(\la w\ra).
\end{equation}

Our next step is to show that if $\mathfrak O$ is an odometer transformation then $\mathfrak O$ has a presentation as an odometer based system with the small word property for some sequence $\la \delta_n:n\in\nn\ra$ tending to  $0$. We do this  by modifying Example \ref{odosolo}.

\begin{lemma}\label{small fingers}
Let $\mathfrak O=\prod_{n\in\nn}\poZ/\poZ_{k_n}$ be an odometer system with invariant measure $\mu$. Then $\mathfrak O$ is isomorphic to  $(\bk, \mu)$ where $\bk$ is the limit of an odometer based construction sequence $\la \mcw_n:n\in\nn\ra$ with $f_n$ tending monotonically to zero exponentially fast; in particular $\sum f_n<\infty$.
\end{lemma}
\pf Let $\mathfrak O$ be an odometer based on $\la k_n:n\in\nn\ra$. Let $n_i$ be a monotone strictly increasing sequence and define $l_i=\prod_{n_{i-1}\le n<n_i}k_n$.  By Proposition \ref{fast growing} $\mathfrak O$ is isomorphic to the odometer based on $\la l_i:i\in \nn\ra$. Thus 
by passing to a subsequence we can assume that:
\[k_{n+1}>3s_n(2^n+1)k_n
.\]
We begin  by letting $\mcw_0=\Sigma=\{a,b,c\}$.\footnote{This construction can be easily modified to work in a 2-letter alphabet, by changing $\mcw_1$ in an \emph{ad hoc} way.}

Suppose that we have constructed $\mcw_n$ and it is enumerated in lexicographical order as $\{w^n_i:1\le i\le s_n\ra$. For each non-identity
 permutation $\sigma$ of 
$\{1, 2, 3, \dots s_n\}$ let $w_\sigma$ be the three-fold concatenation of the words  in $\mcw_n$ in the order given by $\sigma$:
\[w_\sigma=(\prod_{i=1}^{s_n}w^n_{\sigma(i)})^3.\]
Write $k_n=s_n(c_n+3)+d_n$ where $c_n\in \nn, d_n<s_n$ and $\vec{t}=(\prod_{i=1}^{s_n}w^n_i)^{c_n}*\prod_{i=1}^{d_n}w^n_i.$
Finally we let 
\[\mcw_{n+1}=\{w_\sigma\cat \vec{t}:\sigma\in s_n!\}.\]
In words: we begin by making $s_n!-1$ prefixes by concatenating the words in $\mcw_n$ in all possible orders. We then use a single, much longer, suffix to complete each word.

Since each prefix is uniquely readable and comes from a non-trivial permutation $\sigma$,  the words in $\mcw_{n+1}$ are uniquely readable. Moreover any two words in $\mcw_n$ occur with approximately the same frequency in each word in $\mcw_{n+1}$. This precision gets better in a summable way as $n$ increases to  $\infty$. The words in $\mcw_{n+1}$ are clearly concatenations of words in $\mcw_n$.

 By assumption on $k_{n+1}$ the prefix makes up less than $2^{-n}$ portion of a word in $\mcw_{n+1}$. Hence if we let 
 $G=\{x\in \mco:$ for large enough $m, x(m)$ is not in the prefix of any $n$-word$\}$, then as in Example \ref{odosolo}, $G$ is a measure one Borel set and the map 
 $\psi:G\stackrel{1-1}{\longrightarrow}\mco$ continuous.

Since each word in $\mcw_n$  occurs very close to the same number of times in each $\mcw_{n+1}$, the densities of occurrences are all very close to $1/s_n$.   Since  $s_n$ grows as an iterated factorial,  $f_n$  go to zero exponentially.
\qed

If we have an odometer based construction sequence $\la \mcw_n:n\in\nn\ra$ with $f_n\le b_n$ for some sequence $\la b_i:i\in\nn\ra$ going to zero and $\la \delta_i:i\in\nn\ra$ is a sequence of positive numbers less than one, there is a subsequence $\mcv_i=\mcw_{n_i}$ such that $\la \mcv_i:i\in \nn\ra$ has the small word property with respect to $\la \delta_i:i\in\nn\ra$. This subsequence can be chosen continuously in the parameters $\la b_i, \delta_i\ra$. Furthermore,  a tail of any sufficiently fast growing subsequence has the small word property with respect to $\la \delta_n:n\in\nn\ra$. We elaborate on this in the next section.

We now note the following:

\begin{lemma}\label{useful balloon}Let $\mathfrak O$ be an odometer system. Let $\la \mcw_n^{\mathfrak O}:n\in\nn\ra$ be a construction sequence for 
$\mathfrak O$ that has the small word property for $\la \delta_n:n\in\nn\ra$. 
\begin{itemize}
\item If $T:(X,\mu)\to (X,\mu)$ is an ergodic transformation with finite entropy having $\mathfrak O$ as a factor, and $\la \mcw_n^X:n\in\nn\ra$ is the presentation of $X$ as a limit of the odometer based system $\la \mcw^X_n:n\in\nn\ra$ constructed as Theorem \ref{first presentation} as modified in Remark \ref{shifting sands} , then $\la \mcw^X_n:n\in\nn\ra$ has the small word property for $\la \delta_n:n\in\nn\ra$.
\item If $x$ is a Toeplitz sequence with underlying odometer $\mathfrak O$, then the presentation  of the orbit closure $\bl$ of $x$ as  the limit $\bl^*$ of an odometer based construction sequence given in Corollary \ref{much ado} has  the small word property with parameters $\la \delta_n:n\in\nn\ra$.
\end{itemize}
\end{lemma}
\pf In both cases the words in the respective construction sequences were of the form $(u,v)$ where $v$ is in the construction sequence for a presentation of 
$\mathfrak O$.  
\qed

Lemma \ref{useful balloon} reduces the problem of finding presentations of odometer based systems with the small word property to the problem of finding a presentation of the underlying odometer with the small word property. By Lemma \ref{small fingers}, we can do this for a single sequence $\la f_n\ra$ tending to zero. 

\paragraph{The small word property can be arranged continuously} Fix an odometer construction sequence $\la \mcw_n:n\in \nn\ra$, let $n_0=0$ and consider the following game $\mathfrak G(\la \mcw_n:n\in\nn\ra)$.  Let $\la b_n:n\in\nn\ra$ be a sequence with $b_n>f_n$ for all $n$. At round $k\ge 0$:
\begin{itemize}
\item Player I plays $\epsilon_k>0$
\item Player II plays $n_{k+1}>n_k$.
\end{itemize}
Player II wins $\mathfrak G(\la \mcw_n:n\in\nn\ra)$ if and only if $b_{n_{k+1}}<\epsilon_k$ for all $k$.
\bigskip

If is clear that if $b_n$ converges to  $0$ then player II has a winning strategy in $\mathfrak G(\la \mcw_n:n\in\nn\ra)$.  Moreover by Lemma \ref{useful balloon} if $\mathcal S$ is this strategy for  an odometer based presentation 
$\la\mcw_n^{\mathfrak O}:n\in\nn\ra$, then $\mathcal S$ is also a winning strategy for all odometer based presentations   
$\la \mcw_n^X:n\in\nn\ra$ built over $\la\mcw_n^{\mathfrak O}:n\in\nn\ra$. 

In particular we can choose the subsequence $n_k$ Lipshitz continuously in the $\epsilon_k$.

\bibliography{odo_citations}
\bibliographystyle{amsplain}
  
\noindent The first author would like to acknowledge partial support from NSF grant DMS-1700143
      

 \end{document}